\documentclass[a4paper,10pt]{article}
\usepackage{amsmath,amsthm,amsfonts}

\date{}
\RequirePackage[dvips]{graphicx} \textheight 22.5cm
\begin{document}
\centerline{}
\centerline{\bf SOME TYPES OF WEAKLY RICCI SYMMETRIC RIEMANNIAN MANIFOLDS}
\centerline{\bf { Payel Karmakar$^{1}$, Arindam Bhattacharyya$^{2}$ }}
\centerline{$^{1,2}$Department of Mathematics,}
\centerline{Jadavpur University, Kolkata-700032, India.}
\centerline{E-mail: $^{1}$payelkarmakar632@gmail.com, $^{2}$bhattachar1968@yahoo.co.in.}
\newtheorem{Theorem}{\quad Theorem}[section]
\newtheorem{Definition}[Theorem]{\quad Definition}
\newtheorem{Corollary}[Theorem]{\quad Corollary}
\newtheorem{Lemma}[Theorem]{\quad Lemma}
\newtheorem{Example}[Theorem]{\emph{Example}}
\newtheorem{Proposition}[Theorem]{Proposition}
\begin{abstract}
In this paper we discuss when a quasi-conformally flat weakly Ricci symmetric manifold (of dimension greater than 3) becomes a manifold of hyper quasi-constant curvature, a quasi-Einstein manifold and a manifold of quasi-constant curvature. Also we discuss when a pseudo projectively flat weakly Ricci symmetric manifold (of dimension greater than 3) becomes pseudo-quasi constant curvature and a quasi-Einstein manifold, and when a $W_2$-flat weakly Ricci symmetric manifold (of dimension greater than 3) becomes a quasi-Einstein manifold.
\end{abstract}
\section{Introduction}
\
\par{L.Tamassy and T.Q.Binh[11] introduced weakly symmetric Riemannian manifold. A non-flat Riemannian manifold $(M^n,g)(n>2)$ is called weakly symmetric if the curvature tensor $\bar{R}$ of type (0,4) satisfies the condition}\\\
$(\nabla_X \bar{R})(Y,Z,U,V)=A(X)\bar{R}(Y,Z,U,V)+B(Y)\bar{R}(X,Z,U,V)+C(Z)\bar{R}(Y,X,U,V)\\+D(U)\bar{R}(Y,Z,X,V)+E(V)\bar{R}(Y,Z,U,X),~~~~~~~~~~~~~~~~~~~~~~~~~~~~~~~~~~~~~~~~~~~~~~~~~~~~~~~~~~~~~~~~~~~~~~~~~~~~~~~$(1.1)\\\
$\forall$vector fields $X,Y,Z,U,V\in \chi(M^n)$, where $A,B,C,D$ and $E$ are 1-forms (non-zero simultaneously) and $\nabla$ is the operator of covariant differentiation with respect to $g$. The 1-forms are called the associated 1-forms of the manifold and an n-dimensional manifold of this kind is denoted by $(WS)_n$. If the associated 1-forms satisfy $B=C$ and $D=E$, the defining condition of a $(WS)_n$ reduces to the following form[4]\\\
$(\nabla_X\bar{R})(Y,Z,U,V)=A(X)\bar{R}(Y,Z,U,V)+B(Y)\bar{R}(X,Z,U,V)+B(Z)\bar{R}(Y,X,U,V)\\+D(U)\bar{R}(Y,Z,X,V)+D(V)\bar{R}(Y,Z,U,X).~~~~~~~~~~~~~~~~~~~~~~~~~~~~~~~~~~~~~~~~~~~~~~~~~~~~~~~~~~~~~~~~~~~~~~~~~~~~~~$(1.2)\\\
\par{Let $\{e_i\}~~i=1,2,...,n$ be an orthonormal basis of the tangent spaces in a neighbourhood of a point of the manifold. Then setting $Y=V=e_i$ in (1.7) and taking summation over $i,~~1\leq i\leq n$, we get}\\\
$(\nabla_XS)(Z,U)=A(X)S(Z,U)+B(Z)S(X,U)+D(U)S(X,Z)+B(R(X,Z)U)+D(R(X,U)Z),$\\\
where $S$ is the Ricci tensor of type (0,2).$~~~~~~~~~~~~~~~~~~~~~~~~~~~~~~~~~~~~~~~~~~~~~~~~~~~~~~~~~~~~~~~~~~~~~~$(1.3)\\\
\par{A non-flat Riemannian manifold $(M^n,g)(n>2)$ is called weakly Ricci-symmetric if the Ricci tensor $S$ of type (0,2) satisfies the condition}\\\
$(\nabla_XS)(Z,U)=A(X)S(Z,U)+B(Z)S(X,U)+D(U)S(X,Z)~~~~~~~~~~~~~~~~~~~~~~~~~~~~~~~~~~~~~~~~~$(1.4)\\\
\par{From (1.3) it follows that a $(WS)_n(n>2)$ is weakly Ricci symmetric (briefly $(WRS)_n(n>2))$ if [12]}\\\
$~~~~~~~~~~~~~~~~~~~~~~~~~B(R(X,Z)U)+D(R(X,U)Z)=0,~~\forall X,U,Z\in \chi(M^n).~~~~~~~~~~~~~~~~~~~~~~~~~~~~~~~~~~~~~~~~~~~~$(1.5)\\\
\par{Putting $Z=U=e_i$ in (1.4) we get,}\\\
$~~~~~~~~~~~~~~~~~~~~~~~~~dr(X)=rA(X)+B(QX)+D(QX),~~~~~~~~~~~~~~~~~~~~~~~~~~~~~~~~~~~~~~~~~~~~~~~~~~~~~~$(1.6)\\\
where $r$ is the scalar curvature of the manifold.\\\
\par{If a $(WRS)_n(n>2)$ is of zero scalar curvature, then from (1.6),}\\\
$~~~~~~~~~~~~~~~~~~~~~~~~~B(QX)+D(QX)=0,~~\forall X.~~~~~~~~~~~~~~~~~~~~~~~~~~~~~~~~~~~~~~~~~~~~~~~~~~~~~~~~~~~~~~~~~~$(1.7)\\\
\par{From (1.6) we have if a $(WRS)_n(n>2)$ is of non-zero constant scalar curvature, then the 1-form $A$ can be expressed as}\\\
$~~~~~~~~~~~~~~~~~~~~~~~~~A(X)=-\frac{1}{r}[B(QX)+D(QX)],~~\forall X.~~~~~~~~~~~~~~~~~~~~~~~~~~~~~~~~~~~~~~~~~~~~~~~~~~$(1.8)\\\
\par{From (1.4) we obtain,}\\\
$~~~~~~~~~~~~~~~~~~~~~~~~~T(QX)=rT(X),~~\forall X.~~~~~~~~~~~~~~~~~~~~~~~~~~~~~~~~~~~~~~~~~~~~~~~~~~~~~~~~~~~~~~~~~~~~~~~~~$(1.9)\\\
where the vector field $\rho$ is defined by\\\
$~~~~~~~~~~~~~~~~~~~~~~~~~T(X)=g(X,\rho)=B(X)-D(X),~~\forall X.~~~~~~~~~~~~~~~~~~~~~~~~~~~~~~~~~~~~~~~~~~~~~~~~$(1.10)\\\
\par{Also from (1.6) we obtain,}\\\
$~~~~~~~~~~~~~~~~~~~~~~~~~T(Z)S(X,U)-T(U)S(X,Z)=0~~~~~~~~~~~~~~~~~~~~~~~~~~~~~~~~~~~~~~~~~~~~~~~~~~~~~~~~$(1.11)\\\
$\forall$vector fields $X,Z,U$ and $T$ is a 1-form.\\\
\par{B.Das and A.Bhattacharyya studied some types of weakly symmetric Riemannian manifolds in 2011[3]. In 2000, U.C.De et al. discussed about weakly symmetric and weakly Ricci symmetric K-contact manifolds[5] and in 2012, A.A.Shaikh and H.Kundu discussed about weakly symmetric and weakly Ricci symmetric warped product manifolds[10]. Various types of works on weakly Ricci symmetric manifolds is done by S.Jana, A.Shaikh[7] and U.C.De and G.C.Ghosh[6]. Motivated from their work we have established some new results on weakly Ricci symmetric manifolds. In this paper we study quasi-conformally flat, pseudo projectively flat and $W_2$-flat weakly Ricci symmetric manifolds. Here we prove a quasi-conformally flat $(WRS)_n(n>3)$ of non-zero constant scalar curvature is a manifold of hyper quasi-constant curvature and this manifold of non-vanishing scalar curvature is a quasi-Einstein manifold and a manifold of quasi-constant curvature with respect to the 1-form $T$ defined by $T(X)=B(X)-D(X)\neq0~~\forall X$, where $B,D$ are 1-forms(non-zero simultaneously). Also we obtain that a pseudo projectively flat $(WRS)_n(n>3)$ with non-zero constant scalar curvature is a manifold of pseudo quasi-constant curvature and with non-vanishing scalar curvature is a quasi-Einstein manifold with respect to $T$. In the last section we prove that a $W_2$-flat $(WRS)_n(n>3)$ of non-vanishing scalar curvature is a quasi-Einstein manifold with respect to $T$.}\\\
\par{A quasi-conformal curvature tensor $C^*$ is defined by[14]}\\\
 $C^*(X,Y)Z= aR(X,Y)Z+ b[S(Y,Z)X-S(X,Z)Y+ g(Y,Z)QX - g(X,Z)QY]\\-\frac{\gamma}{n}[\frac{a}{n-1}+2b][g(Y,Z)X-g(X,Z)Y],~~~~~~~~~~~~~~~~~~~~~~~~~~~~~~~~~~~~~~~~~~~~~~~~~~~~~~~~~~~~~~~~~~~~~~~~~~~~~~~~~~~$(1.12)\\\
 where $a,b$ are constants, $g$ is the Riemannian metric and $R,Q,\gamma$ are the Riemannian curvature tensor of type (1,3), the Ricci operator defined by $g(QX,Y)=S(X,Y)$ and the scalar curvature respectively.\\\
 \par{Chen and Yano[2] introduced a Riemannian manifold $(M^n,g)(n>3)$ of quasi-constant curvature which is conformally flat and its curvature tensor $\bar{R}$ of type (0,4) is of the form}\\\
 $\bar{R}(X,Y,Z,W)=a[g(Y,Z)g(X,W)-g(X,Z)g(Y,W)]+b[g(X,W)A(Y)A(Z)-g(X,Z)A(Y)A(W)\\+g(Y,Z)A(X)A(W)-g(Y,W)A(X)A(Z)],~~~~~~~~~~~~~~~~~~~~~~~~~~~~~~~~~~~~~~~~~~~~~~~~~~~~~~~~~~~~~~~~~~~~~~~~~~~~~~~~~~~~~~~~~~~~~~~~~$(1.13)\\\
 where $a,b$ are non-zero scalars.\\\
 \par{A Riemannian manifold $(M^n,g)(n>3)$ is said to be of hyper quasi-constant curvature if it is conformally flat and its curvature tensor $\bar{R}$ of type (0,4) satisfies the condition[9]}\\\
 $\bar{R}(X,Y,Z,W)=a[g(Y,Z)g(X,W)-g(X,Z)g(Y,W)]+g(X,W)P(Y,Z)-g(X,Z)P(Y,W)\\+g(Y,Z)P(X,W)-g(Y,W)P(X,Z),~~~~~~~~~~~~~~~~~~~~~~~~~~~~~~~~~~~~~~~~~~~~~~~~~~~~~~~~~~~~~~~~~~~~~~~~~~~~~~~~~~~$(1.14)\\\
 where $a$ is a non-zero scalar and $P$ is a tensor of type (0,2).\\\
 \par{A pseudo projective curvature tensor $\bar{P}$ is defined by[13]}\\\
 $\bar{P}(X,Y)Z=aR(X,Y)Z+b[S(Y,Z)X-S(X,Z)Y]-\frac{r}{n}[\frac{a}{n-1}+b][g(Y,Z)X-g(X,Z)Y],~~~~~~~~~~~~$(1.15)\\\
 where $a,b$ are constants such that $a,b\neq 0$; $R,S,r$ are the Riemannian curvature tensor, the Ricci tensor and scalar curvature respectively.\\\
 \par{A Riemannian manifold $(M^n,g)(n>3)$ is said to be of pseudo quasi-constant curvature if it is pseudo projectively flat and its curvature tensor $\bar{R}$ of type (0,4) satisfies[3]}\\\                                                                                                                                                                                                                                                                                                                                                                              $\bar{R}(X,Y,Z,W)=a[g(Y,Z)g(X,W)-g(X,Z)g(Y,W)]+P(Y,Z)g(X,W)-P(X,Z)g(Y,W),~~~~~$(1.16)\\\
 where $a$ is a constant and $P$ is a tensor of type (0,2).\\\
 \par{Pokhariyal and Mishra have defined the $W_2$-curvature tensor on a differential manifold of dimension $n$ is given by[8]}\\\
 $~~~~~~~~~~~~~~~~~~~~~W_2(X,Y)Z=R(X,Y)Z+\frac{1}{n-1}\{g(X,Z)QY-g(Y,Z)QX\}.~~~~~~~~~~~~~~~~~~~~~~~~~~~~~~~~~~~~~~~~~~~~~~~~~$(1.17)\\\
 \par{All these notions will be required in the next sections.}\\\
\
\section{Quasi-conformally flat $(WRS)_n(n>3)$}
\
In this section we prove that a quasi-conformally flat weakly Ricci symmetric manifold $(WRS)_n\\(n>3)$ of non-zero constant scalar curvature is a manifold of hyper quasi constant curvature and this manifold of non-vanishing scalar curvature is a quasi-Einstein manifold and a manifold of quasi-constant curvature with respect to the 1-form $T$ defined by $T(X)=B(X)-D(X)\neq0~~\forall X$, where $B,D$ are 1-forms(non-zero simultaneously).\\\
\
\par{Let $(M^n,g)(n>3)$ be a quasi-conformally flat $(WRS)_n$. Then from (1.12) we obtain,}\\\
$\bar{R}(X,Y,Z,U)=-\frac{b}{a}[S(Y,Z)g(X,U)-S(X,Z)g(Y,U)+g(Y,Z)S(X,U)-g(X,Z)S(Y,U)]\\+\frac{\gamma}{n}[\frac{1}{n-1}+\frac{2b}{a}][g(Y,Z)g(X,U)-g(X,Z)g(Y,U)],~~~~~~~~~~~~~~~~~~~~~~~~~~~~~~~~~~~~~~~~~~~~~~~~~~~~~~~~~~~~$(2.1)\\\
where $g(R(X,Y)Z,U)=\bar{R}(X,Y,Z,U)$.\\\
\par{Putting $X=U=e_i$ in (2.1) where $\{e_i\}$ is an orthonormal basis of the tangent space at each point of the manifold and taking the summation over $i$, where $1\leq i\leq n$, we get}\\\
$~~~~~~~~~~~~~~~~~~~~~~~~~~~~~~~~~~S(Y,Z)=\alpha g(Y,Z),~~~~~~~~~~~~~~~~~~~~~~~~~~~~~~~~~~~~~~~~~~~~~~~~~~~~~~~~~~~~~~~~~~~~~~~~~~~~~~~~~~~~$(2.2)\\\
where $\alpha =\frac{\gamma}{1+\frac{b}{a}(n-2)}[-\frac{b}{a} +\frac{1}{n}(1+\frac{2b(n-1)}{a})].$
\par{From (2.2) we get}\\\
$~~~~~~~~~~~~~~~~~~~~~~~~~~~~~~(\nabla_XS)(Y,Z)=\alpha^\prime d\gamma(X)g(Y,Z),~~~~~~~~~~~~~~~~~~~~~~~~~~~~~~~~~~~~~~~~~~~~~~~~~~~~~~~~~~$(2.3)\\\
where $\alpha^\prime=\frac{1}{1+\frac{b}{a}(n-2)}[-\frac{b}{a} +\frac{1}{n}(1+\frac{2b(n-1)}{a})].$\\\
\par{Using (2.3) and (1.4) we get,}\\\
$(\nabla_XS)(Y,Z)-(\nabla_ZS)(Y,X)\\=S(Y,Z)[A(X)-D(X)]-S(X,Y)[A(Z)-D(Z)]\\=\alpha^\prime[d\gamma(X)g(Y,Z)-d\gamma(Z)g(Y,X)].~~~~~~~~~~~~~~~~~~~~~~~~~~~~~~~~~~~~~~~~~~~~~~~~~~~~~~~~~~~~~~~~~~~~~~~~~~~~~~~$(2.4)\\\
\par{Let $\rho_1,\rho_2,\rho_3$ be the associated vector fields corresponding to the 1-forms $A,B,D$ respectively, i.e.,}
$g(X,\rho_1)=A(X), g(X,\rho_2)=B(X), g(X,\rho_3)=D(X).$\\\
\par{Putting $Z=\rho_2$ in (2.4) we get,}\\\
$[A(X)-D(X)]B(QY)-S(X,Y)[A(\rho_2)-D(\rho_2)]\\=\alpha^\prime[d\gamma(X)B(Y)-d\gamma(\rho_2)g(Y,X)]\\=\alpha^\prime[B(Y)\{\gamma A(X)+B(QX)+D(QX)\}-g(X,Y)\{\gamma A(\rho_2)+B(Q\rho_2)+D(Q\rho_2)\}].~~~~~~~~~~~~~~~~~~~~~~~~~~~$(2.5)\\\
\par{If the manifold has non-zero constant scalar curvature, then by virtue of (1.8), (2.5) yields that}\\\
$~~~~~~~~~~~~~~~~~[A(X)-D(X)]B(QY)-S(X,Y)[A(\rho_2)-D(\rho_2)]=0~~~~~~~~~~~~~~~~~~~~~~~~~~~~~~~~~~~~~~~~~~~~~~~~~~~~~~~~~~~~~~~~~~~$(2.6)\\\
$~~~~~~~~~~~~~\Rightarrow S(X,Y)=\alpha_1A(X)\bar{B}(Y)+\alpha_2D(X)\bar{B}(Y),~~~~~~~~~~~~~~~~~~~~~~~~~~~~~~~~~~~~~~~~~~~~~~~~~~~~~~~~~~~~~~~~~~~~~~~~~~~~~~~~~~$(2.7)\\\
where $\bar{B}(Y)=B(QY)$ and $\alpha_1=\frac{1}{A(\rho_2)-D(\rho_2)}, \alpha_2=-\frac{1}{A(\rho_2)-D(\rho_2)}.$\\\
\par{This leads to the following theorem:}\\\
\\\
\textit{\textbf{Theorem 2.1:}} \textit{In a quasi-conformally flat $(WRS)_n(n>3)$ of non-zero constant scalar curvature, the Ricci tensor $S$ has the form (2.7).}\\\
\\\
\par{Again using (1.8) in (2.7) we get,}\\\
$~~~~~~~~~~~~~~~~~~~~~~~~~~~~~~~~~~~~S(X,Z)=\alpha_1\bar{B}(Z).-\frac{1}{\gamma}[\bar{B}(X)+\bar{D}(X)]+\alpha_2D(X)\bar{B}(Z),~~~~~~~~~~~~~~~~~~~~~~~~~~~~~~~~~~~~~~~~~$(2.8)\\\
where $\bar{D}(X)=D(QX).$
\par{Using (2.8) in (2.1) we obtain,}\\\
$\bar{R}(X,Y,Z,W)=\alpha[g(Y,Z)g(X,W)-g(X,Z)g(Y,W)]+g(X,W)P(Y,Z)-g(X,Z)P(Y,W)\\+g(Y,Z)P(X,W)-g(Y,W)P(X,Z),~~~~~~~~~~~~~~~~~~~~~~~~~~~~~~~~~~~~~~~~~~~~~~~~~~~~~~~~~~~~~~~~~~~~~~~~~~~~~~$(2.9)\\\
where $\alpha$ is a non-zero scalar and $P(Y,Z)=-\frac{b}{a}[-\frac{\alpha_1}{\gamma}\bar{B}(Y)\bar{B}(Z)+\alpha_2D(Y)\bar{B}(Z)-\frac{\alpha_1}{\gamma}\bar{B}(Z)\bar{D}(Y)].$
\par{Hence comparing (2.9) with (1.14) we can state the following theorem:}\\\
\\\
\textit{\textbf{Theorem 2.2:}} \textit{A quasi-conformally flat $(WRS)_n(n>3)$ of non-zero constant scalar curvature is a manifold of hyper quasi-constant curvature.}\\\
\\\
\par{Now putting $U=\rho$ in (1.11) and then using (1.19) we get,}\\\
$T(Z)S(X,\rho)-T(\rho)S(X,Z)=0$\\\
$\Rightarrow T(Z)T(QX)-T(\rho)S(X,Z)=0$\\\
$\Rightarrow \gamma T(Z)T(X)-T(\rho)S(X,Z)=0.~~~~~~~~~~~~~~~~~~~~~~~~~~~~~~~~~~~~~~~~~~~~~~~~~~~~~~~~~~~~~~~~~~~~~~~~~~~~~~~~~~$(2.10)\\\
\par{Let us suppose that a $(WRS)_n(n>3)$ is quasi-conformally flat and of non-zero scalar curvature.}\\\
\par{From (2.10) we have,}\\\
$~~~~~~~~~~~~~~~~~~~~~~~~~~~~~~~~~~~~~~~~~~~~~~~~S(X,Z)=\tilde{\alpha}T(X)T(Z),~~~~~~~~~~~~~~~~~~~~~~~~~~~~~~~~~~~~~~~~~~~~~~~~~~~~~~~~~~~~~~~~~~~~~~~~~~~~~~$(2.11)\\\
where $\tilde{\alpha}=\frac{\gamma}{T(\rho)}.$\\\
\par{Again, a Riemannian manifold is called quasi-Einstein if its Ricci tensor is of the form[1]$-$}\\\
$~~~~~~~~~~~~~~~~~~~~~~~~~~~~~~~~~~~~~~~~~~~~~~~~S=pg+q\omega\otimes \omega,~~~~~~~~~~~~~~~~~~~~~~~~~~~~~~~~~~~~~~~~~~~~~~~~~~~~~~~~~~~~~~~~~~~~~~$(2.12)\\\
where $p,q$ are scalars of which $q\neq 0$ and $\omega$ is a 1-form.\\\
\par{Comparing (2.11) and (2.12) we can state the following theorem:}\\\
\\\
\textit{\textbf{Theorem 2.3:}} \textit{A quasi-conformally flat $(WRS)_n(n>3)$ of non-vanishing scalar curvature is a quasi-Einstein manifold with respect to the 1-form $T$ defined by $T(X)=B(X)-D(X)\neq 0~~\forall X$.}\\\
\\\
\par{Again using (2.11) in (2.1) we obtain,}\\\
$\bar{R}(X,Y,Z,W)=l[g(Y,Z)g(X,W)-g(X,Z)g(Y,W)]+\delta[g(X,W)T(Y)T(Z)-g(X,Z)T(Y)T(W)\\+g(Y,Z)T(X)T(W)-g(Y,W)T(X)T(Z)],~~~~~~~~~~~~~~~~~~~~~~~~~~~~~~~~~~~~~~~~~~~~~~~~~~~~~~~~~~~~~~~~~~~~~~~~~~~~$(2.13)\\\
where $l=\frac{\gamma}{n}[\frac{1}{n-1}+\frac{2b}{a}]$ and $\delta=-\frac{b}{a}\tilde{\alpha.}$
\par{Hence comparing (2.13) with (1.13) we have the following theorem:}\\\
\\\
\textit{\textbf{Theorem 2.4:}} \textit{A quasi-conformally flat $(WRS)_n(n>3)$ of non-vanishing scalar curvature is a manifold of quasi-constant curvature with respect to the 1-form $T$ defined by $T(X)=B(X)-D(X)\neq 0~~\forall X$.}\\\
\\\
\par{Using the expression of $T$ in (2.13) we have,}\\\
$\bar{R}(X,Y,Z,W)=l[g(Y,Z)g(X,W)-g(X,Z)g(Y,W)]+g(X,Z)\{\delta BD\}(Y,Z)-g(X,Z)\{\delta BD\}(Y,W)\\+g(Y,Z)\{\delta BD\}(X,W)-g(Y,W)\{\delta BD\}(X,Z),$\\\
where $\{\delta BD\}=\delta(BB-BD-DB+DD).$\\\
\par{Comparing the above relation with (1.14) we can state:}\\\
\\\
\textit{\textbf{Corollary 2.1:}} \textit{A quasi-conformally flat $(WRS)_n(n>3)$ of non-zero scalar curvature is a manifold of hyper quasi-constant curvature.}\\\
\\\
\section{Pseudo-projectively flat $(WRS)_n(n>3)$}
\
In this section we obtain that a pseudo projectively flat $(WRS)_n(n>3)$ with non-zero constant scalar curvature is a manifold of pseudo-quasi constant curvature and with non-vanishing scalar curvature is a quasi-Einstein manifold with respect to the 1-form $T$ defined by $T(X)=B(X)-D(X)\neq0~~\forall X$, where $B,D$ are 1-forms(non-zero simultaneously).\\\
\
\par{Let $(M^n,g)(n>3)$ be a pesudo-projectively flat $(WRS)_n$. Then from (1.15) we obtain,}\\\
$\bar{R}(X,Y,Z,W)=-\frac{b}{a}[S(Y,Z)g(X,U)-S(X,Z)g(Y,U)]+\frac{r}{an}[\frac{a}{n-1}+b][g(Y,Z)g(X,U)\\-g(X,Z)g(Y,U)].~~~~~~~~~~~~~~~~~~~~~~~~~~~~~~~~~~~~~~~~~~~~~~~~~~~~~~~~~~~~~~~~~~~~~~~~~~~~~~~~~~~~~~~~~~~~~~~~~~~$(3.1)\\\
\par{Putting $X=U=e_i$ in (3.1) where $\{e_i\}$ is an orthonormal basis of the tangent space at each point of the manifold and taking summation over $i$, $1\leq i\leq n$, we get,}\\\
$~~~~~~~~~~~~~~~~~~~~~~~~~~~~~S(Y,Z)=\alpha g(Y,Z),~~~~~~~~~~~~~~~~~~~~~~~~~~~~~~~~~~~~~~~~~~~~~~~~~~~~~~~~~~~~~~~~~~~~~~~~~~~~~~~~~$(3.2)\\\
where $\alpha =\frac{r}{n}.$\\\
\par{(3.2) indicates that a pseudo-projectively flat manifold is an Einstein manifold.}\\\
\par{Now from (3.2) we get,}\\\
$~~~~~~~~~~~~~~~~~~~~~~~~~~~~~(\nabla_XS)(Y,Z)=\alpha_1dr(X)g(Y,Z),~~~~~~~~~~~~~~~~~~~~~~~~~~~~~~~~~~~~~~~~~~~~~~~~~~~~~$(3.3)\\\
where $\alpha_1=\frac{1}{n}.$\\\
\par{Similarly $(\nabla_ZS)(Y,X)=\alpha_1dr(Z)g(Y,X).~~~~~~~~~~~~~~~~~~~~~~~~~~~~~~~~~~~~~~~~~~~~~~~~~~~~~~~~~~~~~~~~$(3.4)}\\\
\par{Subtracting (3.4) from (3.3) we have,}\\\
$~~~~~~~~~~~~~~~(\nabla_XS)(Y,Z)-(\nabla_ZS)(Y,X)=\alpha_1[dr(X)g(Y,Z)-dr(Z)g(Y,X)].~~~~~~~~~~~~~~~~~~~~~~~~~~~~~~~~~~~~~~~~$(3.5)\\\
\par{Interchanging $X,U$ in (1.4) and then subtracting the resultant from (1.4) we obtain by virtue of (3.5) we obtain,}\\\
$[A(X)-D(X)]S(U,Z)-[A(U)-D(U)]S(X,Z)=\alpha_1[dr(X)g(Z,U)-dr(U)g(Z,X)].~~~~~~~~~~~~~~~~~~~~~$(3.6)\\\
\par{Let $\rho_1,\rho_2,\rho_3$ be the associated vector fields corresponding to the 1-forms $A,B,D$ respectively,}\\\
i.e., $g(X,\rho_1)=A(X),g(X,\rho_2)=B(X),g(X,\rho_3)=D(X).$\\\
\par{Substituting $U$ by $\rho_2$ in (3.6) and then using (1.6) we get,}\\\
$[A(X)-D(X)]B(QZ)-[A(\rho_2)-D(\rho_2)]S(X,Z)=\alpha_1[B(Z)\{rA(X)+B(QX)+D(QX)\}\\-g(X,Z)\{rA(\rho_2)+B(Q\rho_2)+D(Q\rho_2)\}].~~~~~~~~~~~~~~~~~~~~~~~~~~~~~~~~~~~~~~~~~~~~~~~~~~~~~~~~~~~~~~~~~~~~~~~~~~~~~~~~~~~~~~~$(3.7)\\\
\par{If the manifold has non-zero constant scalar curvature, then by virtue of (1.8), (3.7) yields that}\\\
$[A(\rho_2)-D(\rho_2)]S(X,Z)-[A(X)-D(X)]B(QZ)=0~~~~~~~~~~~~~~~~~~~~~~~~~~~~~~~~~~~~~~~~~~~~~~~~~~~~~~~~~~~~~~~~~~~~~~~~~$(3.8)\\\
$\Rightarrow S(X,Z)=\alpha_2A(X)\bar{B}(Z)+\alpha_3D(X)\bar{B}(Z),~~~~~~~~~~~~~~~~~~~~~~~~~~~~~~~~~~~~~~~~~~~~~~~~~~~~~~~~~~~~~~~~~~~~~~$(3.9)\\\
where $\bar{B}(Z)=B(QZ),\alpha_2=\frac{1}{A(\rho_2)-D(\rho_2)},\alpha_3=-\frac{1}{A(\rho_2)-D(\rho_2)}.$\\\
\par{This leads to the following:}\\\
\\\
\textit{\textbf{Theorem 3.1:}} \textit{In a pseudo-projectively flat $(WRS)_n(n>3)$ of non-zero constant scalar curvature, the Ricci tensor $S$ has the form (3.9).}\\\
\\\
\par{Again using (1.8) in (3.9) we obtain,}\\\
$~~~~~~~~~~~~~~~~~~~~~~~~~~~S(X,Z)=-\frac{\alpha_2}{r}\bar{B}(Z)[\bar{B}(X)+\bar{D}(X)]+\alpha_3D(X)\bar{B}(Z).~~~~~~~~~~~~~~~~~~~~~~~~~~~~~~~~~~~~~~~~~~~$(3.10)\\\
\par{Using (3.10) in (3.1) we have,}\\\
$\bar{R}(X,Y,Z,W)=\alpha[g(Y,Z)g(X,W)-g(X,Z)g(Y,W)]+g(X,W)P(Y,Z)-g(Y,W)P(X,Z),~~~$(3.11)\\\
where $\alpha=\frac{r}{an}[\frac{a}{n-1}+b]$ and $P(Y,Z)=-\frac{b}{a}[-\frac{\alpha_2}{r}\bar{B}(Z)\{\bar{B}(Y)+\bar{D}(Y)\}+\alpha_3D(Y)\bar{B}(Z)].$\\\
\par{Comparing (3.11) with (1.16) we can state the following theorem:}\\\
\\\
\textit{\textbf{Theorem 3.2:}} \textit{A pseudo-projectively flat $(WRS)_n(n>3)$ of non-zero constant scalar curvature is a manifold of pseudo quasi-constant curvature.}\\\
\\\
\par{Now putting $U=\rho$ in (1.11) and then using (1.9) we get,}\\\
$~~~~~~~~~~~~~~~~~~~~~~~~~~rT(X)T(Z)-T(\rho)S(X,Z)=0.~~~~~~~~~~~~~~~~~~~~~~~~~~~~~~~~~~~~~~~~~~~~~~~~~~~~~~~~~~~~~~~~~~~~$(3.12)\\\
\par{From (3.12) we have,}\\\
$~~~~~~~~~~~~~~~~~~~~~~~~~~~~~~~~~~~S(X,Z)=\frac{r}{T(\rho)}T(X)T(Z).~~~~~~~~~~~~~~~~~~~~~~~~~~~~~~~~~~~~~~~~~~~~~~~~~~~~$(3.13)\\\
\par{Now $T(\rho)\neq 0$ for if $T(\rho)=0$, then (3.13) implies that $rT(X)T(Z)=0$ which implies that $r=0$ (as $T(X)\neq0~~\forall X$), which is a contradiction.}\\\
\par{Hence $S(X,Z)=\tilde{\alpha_1}T(X)T(Z),~~~~~~~~~~~~~~~~~~~~~~~~~~~~~~~~~~~~~~~~~~~~~~~~~~~~~~~~~~~~~~~~~~~~~~~~~~~~~~$(3.14)}\\\
where $\tilde{\alpha_1}$ is a non-zero scalar.\\\
\par{Hence comparing (3.14) with (2.12) we can state that:}\\\
\\\
\textit{\textbf{Theorem 3.3:}} \textit{A pseudo-projectively flat $(WRS)_n(n>3)$ of non-vanishing scalar curvature is a quasi-Einstein manifold with respect to the 1-form $T$ defined by $T(X)=B(X)-D(X)\neq 0~~\forall X.$}\\\
\\\
\par{Again using (3.14) in (3.1) we obtain,}\\\
$\bar{R}(X,Y,Z,W)=\gamma_1[g(Y,Z)g(X,W)-g(X,Z)g(Y,W)]+\delta_1[T(Y)T(Z)g(X,W)\\-T(X)T(Z)g(Y,W)],~~~~~~~~~~~~~~~~~~~~~~~~~~~~~~~~~~~~~~~~~~~~~~~~~~~~~~~~~~~~~~~~~~~~~~~~~~~~~~~~~~~~~~~~~~~~~~~~~~~~~~~~~~~~~~~~~~~$(3.15)\\\
where $\gamma_1=\frac{r}{an}[\frac{a}{n-1}+b]$ and $\delta_1=-\frac{b}{a}\tilde{\alpha_1}.$\\\
\\\
\par{Comparing (3.15) with (1.16) we have the following theorem:}\\\
\\\
\textit{\textbf{Theorem 3.4:}} \textit{A pseudo-projectively flat $(WRS)_n(n>3)$ of non-vanishing scalar curvature is a manifold of pseudo quasi-constant curvature with respect to 1-form $T$ defined by $T(X)=B(X)-D(X)\neq 0~~\forall X$.}\\\
\\\
\par{Using the expression of $T$ in (3.15), it can be easily seen that}\\\
$\bar{R}(X,Y,Z,W)=\gamma_1[g(Y,Z)g(X,W)-g(X,Z)g(Y,W)]+\{\delta_1BD\}(Y,Z)g(X,W)\\-\{\delta_1BD\}(X,Z)g(Y,W),$\\\
where $\{\delta_1BD\}=\delta_1(BB-BD-DB+DD)$.\\\
\par{Comparing the above relation with (1.16) we can state:}\\\
\\\
\textit{\textbf{Corollary 3.1:}} \textit{A pseudo-projectively flat $(WRS)_n(n>3)$ of non-zero scalar curvature is a manifold of pseudo quasi-constant curvature.}\\\
\\\
\section{$W_2$-flat $(WRS)_n(n>3)$}
\
In this section we prove that a $W_2$-flat $(WRS)_n(n>3)$ of non-vanishing scalar curvature is a quasi-Einstein manifold with respect to the 1-form $T$ defined by $T(X)=B(X)-D(X)\neq0~~\forall X$, where $B,D$ are 1-forms(non-zero simultaneously).\\\
\
\par{Let $(M^n,g)(n>3)$ be a $W_2$-flat $(WRS)_n$. Then from (1.17) we have,}\\\
$~~~~~~~~~~~~~~~~~~~~~~~~~~~~~~~~~~~~~~\bar{R}(X,Y,Z,W)=\frac{1}{n-1}[g(Y,Z)S(X,W)-g(X,Z)S(Y,W)].~~~~~~~~~~~~~~~~~~~~~~~~~~~~~~~~~~~$(4.1)\\\
\par{Putting $X=W=e_i$ in (4.1), where $\{e_i\}$ is an orthonormal basis of the tangent space at each point of the manifold and taking summation over $i$, where $1\leq i\leq n$, we obtain,}\\\
$~~~~~~~~~~~~~~~~~~~~~~~~~~~~~~~~S(Y,Z)=\alpha g(Y,Z),~~~~~~~~~~~~~~~~~~~~~~~~~~~~~~~~~~~~~~~~~~~~~~~~~~~~~~~~~~~~~~~~~~~~~~~$(4.2)\\\
where $\alpha=\frac{r}{n}.$\\\
\par{From (4.2) we get,}\\\
$~~~~~~~~~~~~~~~~~~~~~~~~~~~~~~~(\nabla_XS)(Y,Z)=\alpha^\prime dr(X)g(Y,Z),~~~~~~~~~~~~~~~~~~~~~~~~~~~~~~~~~~~~~~~~~~~~~~~~~~~~~~~~~~~~~$(4.3)\\\
where $\alpha^\prime=\frac{1}{n}.$\\\
\par{Using (4.3) and (1.4) we get,}\\\
$(\nabla_XS)(Y,Z)-(\nabla_ZS)(Y,X)$\\\
$=S(Y,Z)[A(X)-D(X)]-S(X,Y)[A(Z)-D(Z)]$\\\
$=\alpha^\prime [dr(X)g(Y,Z)-dr(Z)g(Y,X)].~~~~~~~~~~~~~~~~~~~~~~~~~~~~~~~~~~~~~~~~~~~~~~~~~~~~~~~~~~~~~~~~~~~~~~~~~~~~~~~~~~$(4.4)\\\
\par{Let $\rho_1,\rho_2,\rho_3$ be the associated vector fields corresponding to the 1-forms $A,B,D$ respectively,}\\\
i.e., $g(X,\rho_1)=A(X),g(X,\rho_2)=B(X),g(X,\rho_3)=D(X).$\\\
\par{Putting $Z=\rho_2$ in (4.4) we get,}\\\
$[A(X)-D(X)]B(QY)-[a(\rho_2)-D(\rho_2)]S(X,Y)$\\\
$=\alpha^\prime [dr(X)B(Y)-dr(\rho_2)g(Y,X)]$\\\
$=\alpha^\prime[B(Y)\{rA(X)+B(QX)+D(QX)\}-g(X,Y)\{rA(\rho_2)+B(Q\rho_2)+D(Q\rho_2)\}].~~~~~~~~~~~~~~~~~~~~~~$(4.5)\\\
\par{If the manifold has non-zero constant scalar curvature, then by virtue of (1.8), (4.5) yields that}\\\
$[A(X)-D(X)]B(QY)-[A(\rho_2)-D(\rho_2)]S(X,Y)=0$\\\
$\Rightarrow S(X,Y)=\alpha_1A(X)\bar{B}(Y)+\alpha_2D(X)\bar{B}(Y),~~~~~~~~~~~~~~~~~~~~~~~~~~~~~~~~~~~~~~~~~~~~~~~~~~~~~~~~~~~~~~~~~~~~~~~~~~~~~~~~~~~~~~~~$(4.6)\\\
where $\bar{B}(Y)=B(QY),\alpha_1=\frac{1}{A(\rho_2)-D(\rho_2)},\alpha_2=-\frac{1}{A(\rho_2)-D(\rho_2)}.$\\\
\par{This leads to the following:}\\\
\\\
\textit{\textbf{Theorem 4.1:} In a $W_2$-flat $(WRS)_n(n>3)$ of non-zero constant scalar curvature, the Ricci tensor $S$ has the form (4.6).}\\\
\\\
\par{Putting $U=\rho$ in (1.11) and then using (1.9) we get,}\\\
$~~~~~~~~~~~~~~~~~~~~~~~~rT(X)T(Z)=T(\rho)S(X,Z).~~~~~~~~~~~~~~~~~~~~~~~~~~~~~~~~~~~~~~~~~~~~~~~~~~~~~~~~~~~~~~~~~~~~~~~~~~~~~~~~~~~~$(4.7)\\\
\par{Let us now suppose that a $(WRS)_n(n>3)$ is $W_2$-flat and of non-zero scalar curvature.}\\\
\par{From (4.7) we have,}\\\
$~~~~~~~~~~~~~~~~~~~~~~~S(X,Z)=\frac{r}{T(\rho)}T(X)T(Z).~~~~~~~~~~~~~~~~~~~~~~~~~~~~~~~~~~~~~~~~~~~~~~~~~~~~~~~~~~~~~~~~~~~~~~~~~~~~~~~~$(4.8)\\\
\par{Now $T(\rho)\neq 0$ as if $T(\rho)=0$ then (4.8) implies that $rT(X)T(Z)=0\Rightarrow r=0$ (as $T(X)\neq 0~~\forall X$), which is a contradiction.}\\\
\par{Hence $S(X,Z)=\tilde{\alpha_1}T(X)T(Z),~~~~~~~~~~~~~~~~~~~~~~~~~~~~~~~~~~~~~~~~~~~~~~~~~~~~~~~~~~~~~~~~~~~~~~~~~~~~~~~~~~$(4.9)}\\\
where $\tilde{\alpha_1}$ is a non-zero scalar.\\\
\par{Hence comparing (4.9) with (2.12) we can state that:}\\\
\\\
\textit{\textbf{Theorem 4.2:} A $W_2$-flat $(WRS)_n(n>3)$ of non-vanishing scalar curvature is a quasi-Einstein manifold with respect to the 1-form $T$ defined by $T(X)=B(X)-D(X)\neq0~~\forall X$.}\\\
\\\
\\\
\textbf{Acknowledgement:}  The first author has been sponsored by University Grants Commission(UGC) Junior Research Fellowship, India. \textit{UGC-Ref. No.}: 1139/(CSIR-UGC NET JUNE 2018).\\\
\\\

\end{document}